\documentclass[12pt]{amsart}

\title{Tennenbaum at Penn and Rochester}
\author{Melvyn B. Nathanson}
\address{Department of Mathematics\\
Lehman College (CUNY)\\
Bronx, New York 10468}
\email{melvyn.nathanson@lehman.cuny.edu}

\begin{document}

\begin{abstract}
Talk at the Logic Conference in Memory of Stanley Tennenbaum at the CUNY Graduate Center  on April 7, 2006.
\end{abstract}

\maketitle

I turned 60 last year, and there was a conference in combinatorial and additive number theory here at the Graduate Center in honor of the sad event.  At such times you get nostalgic if not maudlin, and think back about your life and how you got to wherever you are.  I remembered Stanley Tennenbaum,  who had played a decisive role in my life when I was an undergraduate.  I had not seen or even heard anything about him for several decades, but I decided to send him a thank you note.  I asked the logicians at CUNY if they knew where he lived, and learned that I could contact him through his daughter Susan in Rochester.  Googol and anywho.com provided the address, and I send him a note along with the conference poster, and thanked him for his advice and encouragement when I was a senior at Penn, since the fact that I became a mathematician was due entirely to his influence and example.  Between the time I mailed the note and it arrived in Rochester, Stanley died in Princeton.  The organizers of this conference asked if I would say a few words about his influence on me, and I am happy to do so.  Nothing I say has any mathematical, logical, philosophical, or other significance.  These are simply recollections.  They go back 40 years, and they are as I remember them, which means that they are not necessarily true.  I wrote these notes a few days ago, but I am reassured, after hearing today's other speakers, that my memories are consistent with what others have recounted.

I met Tennenbaum in 1964 at the beginning of the fall semester at the University of Pennsylvania, where he was a visiting professor in the Department of Philosophy.   I was an undergraduate major in philosophy.  It was my senior year, and he was teaching the senior seminar in philosophy.  The main text was Plato's \emph{Theatetus}.  My family lived in Philadelphia, but I lived on campus.  For some reason, I did not have an apartment until a few weeks into the semester, and I was commuting from my home  in West Oak Lane.  Stan and his family lived not far away, and some mornings I would go to his house and drive with him to Penn.  His son Jonathan was attending my old public high school, Central High School, which was something like Styvesant or Bronx Science in New York, except that Central was all boys and its female equivalent, Girls High School, was a block away.   I had a lot of friends at Central and at Penn who went on to get Ph.D.s in mathematics, but at the time I had no interest in mathematics.  The high school calculus teacher was Bernie Warshaw.  He was legendary in Philadelphia, and he really taught calculus.  We had to memorize the epsilon-delta definitions for limits, continuity, and so on,  just as in English classes we had to memorize Shakespeare.  Memorization is greatly underrated.   I think that one of the reasons my eventual transition into mathematics was painless was because of Warshaw's class.  Many of his students went on to become serious  mathematicians.   Any boy who grew up in Philadelphia  and became a mathematician had probably gone to Central and studied calculus with Bernie Warshaw.   

Although I was a philosophy major, I was very interested in biology, and had taken all kinds of course in organic chemistry, physical chemistry, biochemistry, molecular genetics, and biophysics.  I had only taken one year of calculus as a freshman, but Stanley decided I needed to know more mathematics.  He advised me to take two mathematics courses while I was a senior.  One was standard linear algebra with the excellent book \emph{Linear Algebra} by Hoffman and Kunze.  It is still in print.  I tried to use it in a course at Lehman a few years ago, but the book is too hard for our students.   The second was a course in differential equations taught by David Shale; the text was Pontryagin's  \emph{Ordinary Differential Equations}.  At one point Stanley actually told me that instead of reading philosophy for the seminar, I should learn the chapter on Fourier series in Courant's \emph{Calculus}.

Somehow I got caught up in the intellectual and social life of the philosophy department at Penn.  The Penn philosophy department was very ``logical," and besides the usual survey courses on the history of philosophy there was a lot of Frege, G{\" o}del, Quine, and the philosophy of science and mathematics.  The ``great man'' in the department was Nelson Goodman, and I took his course in epistemology.  One of the texts was his little book \emph{Fact, Fiction, and Forecast}, which includes a discussion of the color ''grue,'' which is green before time $t$ and blue thereafter.  The question is:  How do you distinguish the color grue from the color green?   A few years later Goodman moved to Harvard.  I tried to explain grue to the preppies who rowed crew with me at Penn, and they thought I was weird, which in retrospect means that I satisfied one  prerequisite to becoming a mathematician.  There was an odd assortment of faculty and graduate students at Penn who hung out together.  The graduate students included Jay Hullett and Robert Schwartz, and the faculty included not just philosophers but also mathematicians like Peter Freyd and engineers like Robert McNaughton, who were interested in logic and computability issues.  

 I remember that one day as we were driving from North Philadelphia to Penn, Stanley told me the proof of the theorem on invariance of dimension in a vector space.  It was a standard exchange argument to show that the number of vectors in a linearly independent set  cannot exceed the number of vectors in a spanning set for the vector space.  I understood the proof, and Stanley declared that I could become a mathematician.  
 
 I never thought at all about careers, nor about what I might do after Penn.  The Jewish favorite is medicine, but my parents never pushed me in that direction.  It was the family business, however, and I had actually taken vastly more courses in biology, chemistry, and physics than were required for medical school.  I also had a part-time job in the research labs of the Eldridge Reeves Johnson Foundation for Biophysics in Penn's medical school, where one of the faculty, Quentin H. Gibson, had  money in his grant to pay me to putter about in his lab and learn about oxygen uptake by hemoglobin and myoglobin, and the rapid reaction kinetics in enzyme-substrate interactions.   I had taken both the  GRE and the MedCAT, the medical college admission test.  Stanley's advice was that medical school was for people who really wanted to go to medical school and treat patients, and I didn't, so I should go to graduate school.  But in what?  One of the great strengths of American higher education is that undergraduates don't learn very much, at least, not very much about their major.  You study a little of everything in college and a bit more of your major, but you can start graduate school without much knowledge of your discipline.  That's why in America you can major in one subject in college and switch to something else in graduate school.  The curriculum in the rest of the world is very different, of course.  Students are expected to get their general education in high school, and immediately specialize in college, but we also don't teach too much in high school and you have to make up that deficiency in college.  I had studied philosophy and was prepared for medical school; the golden mean was biochemistry or biophysics, and I applied to graduate school.  Stanley showed me the beautiful letter of recommendation that he had written for me, and I wondered how I could live up to it.  I got in everywhere and went to Harvard.  

Studying biophysics was great, but running experiments in a lab for hours every day was less appealing.  At Harvard I studied in the undergraduate library at a table next to the shelf that held Bell's {\em Men of Mathematics}.  No one ever signed it out, and I used to read it whenever I needed a break from my work.  Bell romanticized mathematical life, and influenced many young people into mathematical careers.  I decided to study mathematics and applied to Rochester, becaused that was where Stanley had gone after Penn.  

I showed up in Rochester a few days before the start of the semester.  I didn't know anyone except Stan, so I called him and asked if I could stay with him and his family for a few days.  At the time my request and his invitation seemed perfectly normal, but I realize today, with my own wife and children, what an imposition this must have been, especially on his wife Carol. In the \emph{zeitgeist}  of the 60s, of course, crashing with friends was normal.

Since I did not know any mathematics, the first year at Rochester I took the usual series of undergraduate courses:  real analysis from Rudin, complex variables from Churchill, algebra taught with beautiful notes by Warren May, and a course in number theory taught by Tennenbaum.  It may be that if Stanley had taught an undergraduate course in logic that year, then I would have become a logician, but he taught number theory and I became a number theorist.  

The main text was a magical set of lecture notes by Andr{\' e} Weil, 28 pages of type-written, double-spaced mimeographed pages.  Many years later I was friendly with Morris Schreiber at Rockefeller University, another friend of Stan and someone he had introduced me to.  I gave Weil's notes to Moe, and he passed them on to Walter Kaufman-Buhler, then the mathematics editor of Springer in New York, and Springer subsequently published them in a small softcover volume, \emph{Number Theory for Beginners}, which was supplemented by a set of exercises written by Rosenlicht, who was, evidently, Weil's teaching assistant when he taught the course in Chicago.  The way Stan described the origin of the notes of Weil was as follows:  Weil taught an undergraduate algebra (not number theory) course one summer in Chicago, and these notes were intended to indicate what Weil thought an undergraduate course in abstract algebra should be.  It was in opposition to the then-standard algebra text by Birkhoff and Maclane, which Weil (according to Stanley) thought was terrible.  My impression is that Stanley also did not think much of MacLane.  Somewhat ironically, there is a memorial conference in honor of Saunders Maclane taking place today in Chicago.

One day Stanley decided that the coffee cups in the Rochester math lounge (they were probably paper or styrofoam) were insufficiently elegant for a mathematics department, maybe they did not satisfy the Chicago standard, and we needed better china.  He took me in his car to a store in downtown Rochester to buy good coffee cups and saucers.  We bought some enormous quantity of dishes and schlepped them back to the math department in Stan's car.  Actually, I assumed Stan bought the stuff, but he actually charged it to the math department, and this seems to have created some consternation when the bill arrived.

Stanley resigned from Rochester for reasons that were mysterious to me, but still lived in the city.  One year I had a fellowship that enabled me to travel, and I asked him where were the best places to study number theory.  He told me Cambridge (England), Moscow, and Princeton.   In the next five years I managed to spend an academic year at each of them.   I went to Cambridge for 1969-70, where I took Part III courses with Cassels and Baker, and learned that the British put milk in their tea.  At Cambridge I used my training in philosophy for the first and only time.  Matyasevich has just solved Hilbert's tenth problem, and a preprint of his paper reached Cambridge from Novosibirsk.  I knew a little Russian, and Cassels asked if I would lecture on the paper in the number theory seminar.  Matyasevich's paper, of course, is just an exercise in elementary number theory.  To understand how it solves the undecidability problem for diophantine equations, you have to know a little about recursively enumerable and recursive sets and functions.  I had learned this at Penn, and prefaced my seminar talk with a short introduction to the work of Davis, Putnam, and Robinson.  A few years later Baker saw me at a number theory meeting.  ``What you doing here?" he asked.  ``I thought you were a logician."

I studied with Gel'fand in Moscow in 1972-73, and in 1974-75 I  worked at the Institute for Advanced Study as Assistant to Andr{\' e} Weil, back when the Institute still adhered to the European tradition of letting each permanent professor appoint one visiting member as his ``Assistant.''  When I asked Weil what were my responsibilities as his assistant, he replied, ``Nothing, and conversely."  Both Gel'fand and Weil had fierce reputations, but they were both extraordinarily kind to me as a young mathematician.

Stan often visited the Institute.  I never locked my office door, and on a few occasions he would let himself in at night and sleep there.  A few years later Stan was officially appointed as a visiting member (I think at the instigation of Whitney, with whom Stan had common interests in mathematics education), but after a short time he resigned.  The caretaker of the Institute's apartment complex, a Southerner named Charlie Greb who was a permanent fixture at the Institute, was amused.  After years of hanging around without an office or apartment, Stanley was finally made a member, but after a few weeks he quit.  Charlie Greb didn't understand it.

What I know of the University of Chicago in its golden age is from the anecdotes of Stan Tennenbaum and a few of his friends, like Moe Schreiber.  It may be that their nostalgia distorted reality, and that it was not the intellectually obsessed citadel of the mind that they described.  In support of the picture they painted, however, I recall a review I recently read of a book on admissions to highly selective American universities.  Harvard, Princeton, and Yale admit very smart kids, but it is hard to get in if you can't kick a soccer ball or play squash competitively.  Their graduates become President and generally very successful in the other-than-university world.  Chicago, I think, produced less successful graduates.  One Harvard admissions officer, quoted in the book, said that his great fear was that Harvard would become another University of Chicago, that is, a place concerned with the mind and not the market.  Still, it is comforting to me that once America boasted a university where the intellectual mission was central.  Stan Tennenbaum was the embodiment of that message to me when I was a student, and that was, I think, the reason I responded so strongly to him.  I feel enormously fortunate that 40 years ago I took his seminar on \emph{Theatetus}.

\end{document}